\documentclass[12pt]{article}
\usepackage{amsmath,amssymb,amsthm,epsfig}

\title{On a surface formed by randomly gluing  together polygonal discs} 
\author{Sergei Chmutov and Boris Pittel
\hbox{ }\\
\small 
Department of Mathematics, Ohio State University,\\[-0.8ex]
\small 
231 West 18th Avenue, Columbus, OH 43210, USA\\[-0.8ex]
\small 
\texttt{chmutov@math.ohio-state.edu and bgp@math.ohio-state.edu}}

\begin{document}
\maketitle
{\center\small Mathematics Subject Classifications: 05C80, 05C30, 05A16, 05E10, 34E05, 60C05}
{\center\small Keywords: surfaces, polygonal discs, random permutations, irreducible characters, Euler characteristic, genus, limit distributions}
\def\si{\par\smallskip\noindent}
\def\bi{\par\bigskip\noindent}
\def\pr{\text{ P\/}}
\def\ex{\text{E\/}}
\def\de{\delta}
\def\eps{\varepsilon}
\def\la{\lambda}
\def\a{\alpha}
\def\be{\beta}
\def\de{\Delta}
\def\sig{\sigma}
\def\ga{\gamma}
\def\part{\partial}
\def\Cal{\mathcal}
\def\var{\text{Var\/}}

\def\geq{\geqslant}
\newcommand{\ig}{\includegraphics}
\newcommand{\rb}{\raisebox}
\newcommand\risS[6]{\rb{#1pt}[#5pt][#6pt]{\begin{picture}(#4,15)(0,0)
  \put(0,0){\ig[width=#4pt]{#2.eps}} #3
     \end{picture}}}
\newcommand\risSpdf[6]{\rb{#1pt}[#5pt][#6pt]{\begin{picture}(#4,15)(0,0)
  \put(0,0){\ig[width=#4pt]{#2.pdf}} #3
     \end{picture}}}

\newtheorem{Theorem}{Theorem}[section]
\newtheorem{Lemma}{Lemma}[section]
\newtheorem{Proposition}{Proposition}[section]
\newtheorem{Corollary}{Corollary}[section]
\newtheorem{Conjecture}{Conjecture}[section]
\numberwithin{equation}{section}

\begin{abstract}
Starting with a collection of $n$ oriented polygonal discs, with an even number $N$ of sides
in total,  we generate a random oriented surface by randomly matching the sides of discs and 
properly gluing them together. Encoding the surface in a random permutation $\gamma$ of $[N]$, we use the Fourier transform on $S_N$ to show that $\gamma$ is asymptotic to the permutation distributed 
uniformly on the alternating group $A_N$ ($A_N^c$ resp.) if $N-n$ and $N/2$ are of the
same (opposite resp.) parity. We use this to prove a local central limit theorem  for the number of vertices on the surface, whence for its Euler characteristic $\chi$. We also show that with high probability the random surface consists of a single component, and thus has a well-defined
genus  $g=1-\chi/2$, which is asymptotic to a Gaussian random variable, with
mean $(N/2-n-\log N)/2$ and variance $(\log N)/2$.
\end{abstract}

\section{Introduction and main results}

In this paper we study random surfaces obtained by gluing, uniformly at random, sides of $n$ polygons with various (not necessarily equal) number of sides. We call this scheme of generating a  surface the {\it map model}.  (A model dual to the map model is very important for algebraic geometry
\cite{LZ}. It can be generalized to hypermaps;  in \cite {ChVT} it is called $\sigma$-{\it model}.)
In the map model the interiors of polygons represent countries ({\it faces}); the glued sides represent boundaries between countries ({\it edges}). Thus the map model can be considered as a graph embedded into the surface such that the faces correspond to the original polygons.

This model generalizes the random map model of N.~Pippenger and K.~Schleich \cite{PippengerSchleich} where all the polygons are triangles, a model motivated by studies in quantum gravity. In particular,  for the Euler characteristic $\chi$ of the randomly
triangulated  surface they proved that $\ex[\chi]=n/2-\log n+O(1)$, $\var(\chi)=\log n+O(1)$, and made startlingly sharp conjectures regarding the remainder terms $O(1)$, based on simulations  and results for similar models. The case when the number of sides of all polygons are equal, gluings of $k$-gons ($k\geq 3$), was considered by A.~Gamburd in \cite{Gamburd}. His breakthrough result was that  (for $2\text{lcm}\{2,k\}\,|\,3n$),  the underlying random permutation of polygons sides was asymptotically uniform on the alternating subgroup $A_{3n}$, implying, for instance, that $\chi$ was asymptotic, in
distribution, to $n/2$ minus ${\cal N}(\log n,\log n)$, the Gaussian variable, with mean and
variance equal $\log n$. K. Fleming and N. Pippenger \cite{FlemingPippenger} used Gamburd's
result to prove sharp asymptotic formulas for the first four moments of the Euler characteristic
$\chi$, in particular confirming their earlier conjectures for $k=3$. 
Another special case of this model is when there is only one polygon whose sides are glued in pairs. This case is well studied, popular, and important in combinatorics and the theory of moduli spaces of algebraic curves. The classical paper of J.~Harer and D.~Zagier \cite{HarerZagier} 
solved the difficult problem of enumerating the resulting surfaces by genus. Their result was used in \cite{ChmutovPittel} to determine the limiting genus distribution for the surface chosen uniformly
at random from all such surfaces.

The sides of the polygons are glued in pairs. So the total number of sides $N$ of all polygons must be even, and the resulting map will have $N/2$ edges.
We also assume that all polygons are directed and that in each glued pair the edges are
directed opposite-wise. Thus the resulting surface will be oriented.

The map model can described in terms of permutations. Label $e$'s the directed sides (edges) of all polygons by numbers from $[N]:=\{1,2,\dots,N\}$; $e_j$ will denote the edge labeled $j$. Let $n_j$ be the number of polygons with $j$ sides, $j$-gons, and let
$J$ stand for the set of all possible numbers of sides of our $n=\sum_jn_j$ polygons, so that
$\sum_{j\in J}jn_j=N$ and each map will have $n$ faces. 
We define the permutation $\alpha$ of $[N]$ as follows: $\alpha(e_j)=e_k$ if $e_k$ follows,
immediately,  $e_j$
in one of the $n$ directed polygons. Thus $\alpha$ has $n$ cycles, each cycle consisting of the edges of the attendant polygon listed according to the polygon orientation.
The set of all such $\alpha$'s is
the conjugacy class ${\cal C}_J$ of permutations of $[N]$ with $n_j$ cycles of length $j$.  A gluing itself is encoded in the permutation $\beta$ which is a product of transpositions of edges  that are glued to each other; those $\beta$'s are
all $(N-1)!!$ elements of the conjugacy class ${\cal C}_2$ of permutations of $[N]$ with cycles of length $2$ only. 

Here is how a given pair of permutations $\alpha$, $\beta$ induces the corresponding surface.
The first edge $e_1$ is glued to the edge $\beta(e_1)$;  the edge $\beta(e_1)$ is followed by the edge
$e_2=\alpha(\beta(e_1))=(\alpha\beta)(e_1)$ in the directed polygon that contains $\beta(e_1)$. 
Next $e_2$ is glued to $\beta(e_2)$ followed  by $e_3=\alpha(\beta(e_2))=(\alpha\beta)(e_2)$ in the cycle that contains $\beta(e_2)$, and so on, producing a sequence of edges $e_1,e_2,\dots$, 
whose {\it tails\/} are lumped together as a single vertex. Since $\alpha\beta$ is a permutation of $[N]$, the sequence $e_1,e_2,\dots$ eventually loops back on the starting edge $e_1$, forming a cycle $e_1\to e_2\to\cdots\to e_m\to e_1$ of $\alpha\beta$, see the picture:
$$\risSpdf{-75}{cycles}{\put(40,17){\rotatebox{80}{$\scriptstyle e_1$}}
    \put(50,25){\rotatebox{-67}{$\scriptstyle \beta(e_1)$}}
    \put(65,49){\rotatebox{-15}{$\scriptstyle \alpha(\beta(e_1))=:e_2$}}
    \put(78,58){\rotatebox{10}{$\scriptstyle \beta(e_2)$}}
    \put(48,72){\rotatebox{65}{$\scriptstyle \alpha(\beta(e_2))=:e_3$}}
    }{100}{40}{80}\qquad
\risSpdf{-70}{cycles2}{\put(82,20){$\scriptstyle e_1$}\put(115,32){$\scriptstyle e_2$}
         \put(110,64){$\scriptstyle e_3$}}{140}{0}{0}
$$
Likewise, starting from the first edge not in this cycle, i. e. distinct from $e_1,\dots, e_{m}$,
we obtain an {\it independent\/} cycle containing this edge, that determines another vertex of the map.
Proceeding in this fashion, we eventually partition the edge set into disjoint subsets, each 
associated with its own vertex of the map. Clearly, the number of those subsets, i. e. 
the number of vertices $V_N$, equals the number of cycles of $\gamma:=\alpha\beta$. 

Also it is obvious that the connected components of the
resulting surface correspond to the orbits of the subgroup generated by permutations $\alpha$ and $\beta$. For each such orbit we know the number of faces (which is equal to the number of cycles of $\alpha$ restricted to the orbit), the number of edges (which is equal to the  number of cycles of $\beta$ restricted to the orbit), and the number of vertices (which is equal to the  number of cycles of $\gamma$ restricted to the orbit). Thus we know the Euler characteristic of the corresponding connected component. Since the Euler characteristic is a complete topological invariant of connected oriented surfaces, the permutations  $\alpha$ and $\beta$ completely deteremine the topology of the surface.

For example, consider two oriented squares with labeled sides: 
$$\risSpdf{-45}{2-squares}{}{150}{0}{40}\qquad
$$
The following three gluings result in a sphere, a torus, and two tori respectively.

$1)\qquad \risSpdf{-35}{S2}{\put(220,60){$\alpha=(1234)(5678)$}
        \put(220,40){$\beta=(15)(28)(37)(46)$}
        \put(220,20){$\gamma=(16)(25)(38)(47)$}}{200}{45}{40}$

$2)\quad \risSpdf{-80}{T2}{\put(230,90){$\alpha=(1234)(5678)$}
        \put(230,70){$\beta=(15)(24)(37)(68)$}
        \put(230,50){$\gamma=(1652)(3874)$}}{210}{60}{80}$

$3)\quad \risSpdf{-80}{2T2}{\put(230,90){$\alpha=(1234)(5678)$}
        \put(230,70){$\beta=(13)(24)(57)(68)$}
        \put(230,50){$\gamma=(1432)(5876)$}}{210}{60}{90}$

We define a {\it random surface} as a surface obtained by gluing via the permutations $\alpha$ and $\beta$ that are chosen uniformly at random (uar), and independently of each other, from  the conjugacy classes ${\cal C}_{J}$ and  ${\cal C}_{2}$ respectively.

In Section \ref{s:lim-uniform} we show, Theorem \ref{conj2},  that the probability distribution of the permutation $\gamma$ is asymptotically uniform on $A_N$ ($A_N^c$ resp.) if ${\cal C}_2$ and ${\cal C}_J$ are
of the same (opposite resp.) parity. This generalizes  result from \cite[Theorem 4.1]{Gamburd} and improves the guaranteed rate of convergence from $N^{-1/12}$ to $N^{-1}$. By and large,
we follow \cite{Gamburd} to reduce the problem to Fourier-based analysis of the total variation distance between two probability measures on $S_N$, main tool being a fundamental general
bound due to P. Diaconis and M. Shashahani \cite{DiaconisShashahani}.  At the crucial point, when we need to estimate  a character value of an irreducible representation  on a 
{\it general\/} conjugacy class of $S_N$ (rather than the classes ${\cal C}_{\ell}$  in \cite{Gamburd} treated via a bound discovered by S. Fomin and N. Lulov \cite{FominLulov} in $1997$),  we use a bound  proved recently by M. Larsen and A. Shalev \cite{LarsenShalev}. 

In Section \ref{s:num-vert}, as a corollary of Theorem \ref{conj2},  we state  that  the total variation distance
between $V_N=V_{N,{\cal C}_J}$, the number of vertices on the random surface, and the number of cycles $C_N^e$
($C_N^o$ resp.) in the uniformly random even (odd resp.) permutation of $[N]$ is of order $O(N^{-1})$. Our main result, a local central limit theorem (LCLT) for $V_N$, follows then from a 
LCLT for $C_N$, the number of cycles in the permutation  distributed uniformly on $S_N$, due to
V. Kolchin \cite{Kolchin}. The LCLT for the Euler characteristic $\chi_N=\chi_{N,{\cal C}_J}$ of the random surface follows immediately.

In the last Section \ref{s:num-comp} we discuss the distribution of the number of connected components of the surface. Generalizing the result of Pippenger and Schleich for $J=\{3\}$, 
\cite{PippengerSchleich}, we prove in Theorem
\ref{connect} that the resulting surface is connected with probability $1-O(N^{-1})$. Thus,
with high probability, the genus $g_N=g_{N,J}$ of the random surface is well defined, and using
$g_N=1-\chi_{N}/2$ we obtain a LCLT for $g_{N}$. For a very special case of one polygon, $J=\{N\}$,
this proves a slightly weaker version of our earlier result in \cite{ChmutovPittel}.

\section{Limiting uniformity} \label{s:lim-uniform}

Given $N$, let $J=J(N)$ be a subset of $\{3,4,\dots\}$. Let $\{n_j\}$ be such that $\sum_{j\in J}jn_j=N$. Consider the set of all partitions of $[N]$ into $n=\sum_{j\in J}n_j$ disjoint cycles, with
$n_j$ cycles of lengths $j\in J$.
This set can be viewed as the conjugacy class  ${\cal C}_J={\cal C}_{N,J}$ of all permutations $\alpha\in S_N$ with $n_j$ cycles of length $j\in J$. 
Assuming that $N$ is divisible by $2$, let  $ {\cal C}_{ 2}$ denote 
the conjugacy class of $S_N$ consisting of permutations that are products of $N/2$ disjoint
$2$-cycles.

Let $\alpha$ and $\beta$ be chosen independently of each other, uniformly at random
(uar) from  ${\cal C}_{J}$ and  ${\cal C}_{2}$ respectively, and let $\gamma=\alpha\beta$.

Gamburd \cite{Gamburd} had studied a special case when $J$
is a singleton $\{k\}$, ($k\ge 3$), and  $\alpha,\beta$, whence $\gamma$, are {\it all even\/}.  
(A permutation $\sigma$ of $[N]$ is called even if it has an even number of even cycles, or equivalently if $N$ minus  the number of cycles of $\sigma$ is even.)
\begin{Theorem}(Gamburd)\label{Gam1}  Suppose that $N\to\infty$ through values
divisible by $2\,\text{lcm}\{2,k\}$. Let $P_{\gamma}$ be the probability distribution of $\gamma$ and
let $U$ be the uniform probability measure on the alternating subgroup $A_N$ of even permutations. Let $\|P_{\gamma}-U\|=\|P_{\gamma}-U\|
_{\text{TV}}$ denote the total variation distance between $P_{\gamma}$ and $U$. Then
\begin{equation}\label{-1/12}
\|P_{\gamma}-U\|=O\bigl(N^{-1/12}\bigr).
\end{equation}
\end{Theorem}
\noindent As noted in Fleming and Pippenger \cite{FlemingPippenger}, the original condition $\text{lcm}\{2,k\} | N$ in \cite{Gamburd}
does not guarantee that both $\alpha$ and $\beta$ are even, implying evenness of $\gamma$. 
Namely (assuming $\text{lcm}\{2,k\} | N$): (1) $\beta$ is even (odd resp.), if $4 | N$ ($4\!\not |\,N$
resp.); (2)
if $k$ is even and $2k | N$ 
($2k\!\not |\, N$ resp.), then $\alpha$ is even (odd resp.); (3) if $k$ is odd then $\alpha$ is even.
Thus $\alpha,\beta,\gamma$ are all even iff $2\,\text{lcm}\{2,k\}|N$; $\gamma$ itself is even
iff $\alpha$ and $\beta$ are of the same parity, i. e. iff $2k\, |\, N(k-2)$.

Gamburd proved \eqref{-1/12} by using a character-based  bound, due to
Diaconis and Shashahani \cite{DiaconisShashahani},  for the 
{\it total variation distance\/} between two probability measures (one being uniform) on a general
finite group $G$ in the special case when $G$ was the alternating subgroup $A_N$. We found
that Gamburd's argument can be modified to prove a far more general, and stronger, result  by using a {\it variation\/} of the bound  in \cite{DiaconisShashahani} for the group $S_N$ itself, when the ``uniform'' measure is supported either by $A_N$ or its coset $A_N^c$, dependent upon parity of $\gamma$.
\begin{Theorem}\label{conj2} Uniformly over
all the classes $\cal C_{J}$ with $\min J\ge 3$, $\gamma=\alpha \beta$ is asymptotically uniform over $A_N$ 
(over $A_N^c$ resp.) if ${\cal C}_J$, ${\cal C}_2$ are of the same parity (of opposite  parity resp.),
and more precisely
\begin{equation}\label{-1}
\|P_{\gamma}-U_{A_N}\|=O\bigl(N^{-1}\bigr), \quad \left(\|P_{\gamma}-U_{A_N^c}\|=
O\bigl(N^{-1}\bigr)\,{ resp.}\right),
\end{equation}
$U_{A_N}$, $U_{A_N^c}$ being the probability measures uniform on $A_N$ and $A_N^c$
respectively.
\end{Theorem}
\noindent For $|J|=1$, and even $\alpha$, $\beta$, the (first) bound in \eqref{-1} improves the bound \eqref{-1/12}.\\

\begin{proof} Like the proof of Theorem \ref{Gam1} in \cite{Gamburd}, the starting point  is the already mentioned Diaconis-Shashahani's bound. Let $G$ be a finite group and $P$, $U$ be two probability measures
on $G$, $U$ being uniform, i. e. $U(g)=1/|G|$ for every $g\in G$. Then 
\begin{equation}\label{DS}
\|P-U\|^2\le\frac{1}{4}\sum_{\rho\in \widehat G,\,\rho\neq \text{id}}\text{dim}(\rho)\, \text{tr}\bigl(\hat P(\rho)
\hat P(\rho)^*\bigr);
\end{equation}
here $\hat G$ denotes the set of all irreducible representations $\rho$ of $G$, "id" denotes the
trivial representation, $\text{dim}(\rho)$ is the dimension of $\rho$, and $\hat P(\rho)$
is the matrix value of the Fourier transform of $P$ at $\rho$. This bound followed from
Cauchy-Schwartz inequality
\begin{equation}\label{CS}
4\|P-U\|^2\le |G|\sum_{s\in G}|P(s)-U(s)|^2,
\end{equation}
combined with the Plancherel Theorem
\begin{equation}\label{PT}
|G|\sum_{s\in G}|P(s)-U(s)|^2=\sum_{\rho\in\widehat G}d(\rho)\text{tr}\bigl[(\hat P(\rho)-\hat{U}(\rho))
(\hat P(\rho)-\hat{U}(\rho))^*\bigr],
\end{equation}
and the observation that  {\bf (i)\/} $\hat P(\rho)=\hat U(\rho) = 1$ for $\rho=\text{id}$, and {\bf (ii)\/}
$\hat U(\rho)=0$ for $\rho\neq \text{id}$. 

Now, \eqref{CS}-\eqref{PT} hold for any two measures on $G$, whence for two probability
measures $P_H$ and $U_H$ supported by the same subset $H\subseteq G$. In this case, 
the condition {\bf (i)\/} still holds, and  we get
\begin{equation}\label{interm}
\|P_H-U_H\|^2\le \frac{1}{4}\sum_{\rho\neq \text{id}}d(\rho)\text{tr}\bigl[(\hat{P}_H(\rho)-\hat{U}_H(\rho))(\hat{P}_H(\rho)-\hat{U}_H(\rho))^*\bigr].
\end{equation}
For $G=S_N$, the irreducible representations $\rho$ are labeled by $\la$, where each $\la
$ is a partition $\la= (\la_1\ge \la_2\ge\dots)$ of $N$, $\la\vdash N$ in short, and $\text{dim}(\rho^{\la})=
f^{\la}$, given by the hook formula
\[
f^{\la}=\frac{N!}{\prod_{u\in \la}h(u)}.
\]
Furthermore one-row $\la =\langle N\rangle$ is the identifying label of the trivial representation ``id'', and
one-column $\la=\langle 1^N\rangle$ is the label of the second one-dimensional representation
``$\text{sign}$'', with value $1$ on $A_N$ and value $-1$ on $A_N^c$. In our case $H$ is either $A_N$ or $A_N^c$, so $\text{sign}(\sigma)$ is the same, $\text{sign}(H)$, for all permutations $\sigma\in H$. Consequently, for $\rho=\text{sign}$,
\[
\hat{P}_H(\rho)=\sum_{\sigma\in H}\text{sign}(\sigma)P_H(\sigma)=\text{sign}(H)\sum_{\sigma\in H}
P_H(\sigma)=\text{sign}(H),
\]
and likewise $\hat{U}_H(\rho)=\text{sign}(H)$. Therefore
\begin{equation}\label{la=sign}
\hat{P}_H(\rho)-\hat{U}_H(\rho)=0,\quad (\rho=\text{sign}).
\end{equation}
Consider $\la\neq \langle N\rangle,\, \langle1^N\rangle$. If $\la$ is not self-dual, i. e. $\la\neq \la^
\prime$, then the $\rho^{\la}$ restricted to $A_N$ is a nontrivial  irreducible representation
$\rho$ of $A_N$, whence $\sum_{\sigma\in A_N}\rho^{\la}(\sigma)=0$, (Diaconis \cite{Diaconis}, Ch. 2B, Exer. 3). Of course $\sum_{\sigma\in S_N}\rho^{\la}(\sigma)=0$ too, whence we have
\begin{equation}\label{hat U(H)=0,1}
\hat U_H(\rho^{\la})=\frac{1}{|H|} \sum_{\sigma\in H}\rho^{\la}(\sigma)=0, \quad (\la\neq \la^\prime).
\end{equation}
If $\la=\la^\prime$ then $\rho^{\la}$ restricted to $A_N$ is a direct sum of two irreducible representations
each of dimension $f^{\la}/2$, which exceeds $1$ for $N\ge 5$, because $f^{\la}\ge 6$ for the self-dual $\la$ with $|\la|\ge 5$. Therefore again we have: for $N\ge 5$,
\begin{equation}\label{hat U(H)=0,2}
\hat U_H(\rho^{\la})=\frac{1}{|H|} \sum_{\sigma\in H}\rho^{\la}(\sigma)=0, \quad (\la=\la^\prime,\, |\la|\ge 5).
\end{equation}

Putting together \eqref{interm}, \eqref{la=sign}, \eqref{hat U(H)=0,1} and \eqref{hat U(H)=0,2}, we obtain: for $N\ge 5$ and $H=A_N$ or $H=A_N^c$,
\begin{equation}\label{|PH-UH|}
\|P_H-U_H\|^2\le \frac{1}{4}\sum_{\la\neq \langle N\rangle,\, \langle 1^N\rangle}f^{\la}\,
\text{tr}\bigl[\hat P_H(\rho^{\la})\hat P_H(\rho^{\la})^*\bigr].
\end{equation}

Once \eqref{|PH-UH|} is proved, the next step is essentially the same as in Gamburd's 
argument when $J=\{k\}$, ${\cal C}_k$, ${\cal C}_2$ are both even, implying that $H=A_N$. In our case $P_H=P_{\gamma}=U_{{\cal C}_J}\star U_{{\cal C}_2}$,  and so, by multiplicativity
of the Fourier transform for convolutions,
\[
\hat{P}_H(\rho^{\la}) = \hat {U}_{{\cal C}_J}(\rho^{\la}) \cdot \hat {U}_{{\cal C}_2}(\rho^{\la}).
\]
Since $U_{{\cal C}_J}=|{\cal C}_J|^{-1}\, 1_{{\cal C}_J}$, $U_{{\cal C}_2}
=|{\cal C}_2|^{-1}\, 1_{{\cal C}_2}$ are class functions, each supported by a single
conjugacy class,
\[
\hat {U}_{{\cal C}_J}(\rho^{\la})=\frac{\chi^{\la}({\cal C}_J)}{f^{\la}}\,I_{f^{\la}},\quad
\hat {U}_{{\cal C}_2}(\rho^{\la})=\frac{\chi^{\la}({\cal C}_2)}{f^{\la}}\,I_{f^{\la}};
\]
here $\chi^{\la}$ is the character of $\rho^{\la}$. So
\[
\hat {P}_H(\rho^{\la})=\frac{\chi^{\la}({\cal C}_J)\chi^{\la}({\cal C}_2)}{(f^{\la})^2}\,I_{f^{\la}}
\]
and therefore \eqref{|PH-UH|} becomes
\begin{equation}\label{DSG}
\|P_\gamma-U_H\|^2\le \frac{1}{4}\sum_{\la\neq \langle N\rangle,\, \langle 1^N\rangle}\left(\frac{
\chi^{\la}({\cal C}_J)\chi^{\la}({\cal C}_2)}{f^{\la}}\right)^2.
\end{equation}
With $1/2$ instead of $1/4$ and ${\cal C}_k$ instead of the general ${\cal C}_J$, the RHS of \eqref{DSG} is Gamburd's upper bound for his case. To make use of his bound, Gamburd applied the following estimate due to Fomin and Lulov \cite{FominLulov}: for $N=tn$,
\begin{equation}\label{FL} 
|\chi^{\la}({\cal C}_t)| =O\bigl(N^{1/2-1/(2t)}\bigr)(f^{\la})^{1/t},
\end{equation}
uniformly for all $N$ and $\la$. He used \eqref{FL} for for both $t=2$ and $t>2$. For $|J|>1$
a similar bound for $|\chi^{\la}({\cal C}_J)|$ was not available at that time. More
recently Larsen and Shalev \cite{LarsenShalev} proved a remarkable extension of the Fomin-Lulov bound:
given $m$, uniformly for all  permutations $\sigma$ without cycles of length below $m$,
and partitions $\la$,
\begin{equation}\label{L-S}
|\chi^{\la}(\sigma )|\le (f^{\la})^{1/m+o(1)},\quad N\to\infty.
\end{equation}
(For $m=2$, i. e. for fixed-point-free permutations, this is very similar to a bound conjectured
earlier by Fomin and Lulov.) With this bound applied to both $\chi^{\la}({\cal C}_2)$
and $\chi^{\la}({\cal C}_J)$, the remaining proof of Theorem \ref{conj2} largely, but not entirely,  follows the original Gamburd's argument.

Introduce $\Lambda=\{\la\vdash N: \lambda_1\ge N-6\}$ and write
\begin{align*}
\|P_{\gamma}-U\|^2&\le \frac{1}{2}\sum_{\la\in \widehat{S_N}\atop \la\neq \langle N\rangle,\,\langle
1^N \rangle}\left(\frac{\chi^{\la}({\cal C}_J)\chi^{\la}({\cal C}_2)}{f^{\la}}\right)^2\\
&\le\sum_{\la\vdash N\atop \la_1\le N - 7}\left(\frac{\chi^{\la}({\cal C}_J)\chi^{\la}({\cal C}_2)}{f^{\la}}\right)^2+
\sum_{\la\vdash N\atop \la\in \Lambda}\left(\frac{\chi^{\la}({\cal C}_J)\chi^{\la}({\cal C}_2)}{f^{\la}}\right)^2\\
&=: \Sigma_1+\Sigma_2.
\end{align*}
Consider $\Sigma_1$. By \eqref{L-S},
\begin{align*}
\left(\frac{\chi^{\la}({\cal C}_J)\chi^{\la}({\cal C}_2)}{f^{\la}}\right)^2&\le \left(\frac{(f^{\la})^{1/3+1/2+o(1)}}{f^{\la}}\right)^2\\
&=(f^{\la})^{-1/3+o(1)};
\end{align*}
so using Proposition 4.2 (Gamburd),
\begin{equation}\label{Sigma1=o}
\Sigma_1=\left.O\bigl(N^{-7t}\bigr)\right|_{t=1/3-o(1)}=o(N^{-2}).
\end{equation}
To handle $\Sigma_2$, we use the following bounds. If $a>0$ is fixed, then uniformly
for $\la$ such that $\la_1=N-a$, and ${\cal C}_J$,
\begin{equation}\label{nufixed}
\begin{aligned}
&f^{\la}\ge \binom{N-a}{a}\ge \frac{N^{a}}{2a!},\\
&|\chi^{\la}({\cal C}_2)|=O\left(N^{\lfloor a/2\rfloor}\right),\quad 
|\chi^{\la}({\cal C}_J)|=O\left(N^{\lfloor a/3\rfloor}\right).
\end{aligned}
\end{equation}
For the first line bound see \cite{Gamburd}  equation (4.17). Let us prove the second line bounds.
Consider $|\chi^{\la}({\cal C}_J)|$, for example. ${\cal C}_J$ is a set of all permutations 
 $\sigma\in A_N$ whose cycles are of lengths from $J$, with fixed counts of cycles of each admissible
length. Let $\boldsymbol{\alpha}=(\alpha_1,\alpha_2,\dots)$ be an arbitrary composition formed
by cycle lengths of a permutation $\sigma\in {\cal C}_J$. From Murnaghan-Nakayama rule,
Stanley \cite{Stanley} (Section 7.17, Equation (7.75)),
\[
|\chi^{\la}({\cal C}_J)|\le g^{\la}(\boldsymbol{\alpha}),
\]
where  $g^{\la}(\boldsymbol{\alpha})$ is the total number of ways to empty the diagram $\la$ by
by successive deletion of the rim hooks, one hook at a time, of lengths $\alpha_1,\alpha_2,\dots$.
Let us show that  
\[
g^{\la}(\boldsymbol{\alpha})=O\left(N^{\lfloor a/3\rfloor}\right).
\]
Each of the $g^{\la}(\boldsymbol{\alpha})$ ways to empty $\la$ consists of an ordered sequence of hook deletions not touching
any of the first $a$ cells in the first row, concatenated with an ordered sequence of hook 
deletions, the first of which deletes at least the cell $(1,a)$ from among those $a$ cells, with the remaining deletions
taking place entirely in a remaining corner-subdiagram $\mu$, with $|\mu|\le 2(a-1)$.
So the number of ways
to empty the residual diagram $\mu$ is at most some $S_1(a)=O(1)$, as $a$ is fixed. As for the
first batch of hook deletions, they are deletions of {\it horizontal\/} rim hooks from
the first row, possibly interspersed with deletions of rim hooks from the subdiagram $\nu$ formed
by all the other rows of $\la$. Since $|\nu|=a$, the length of the subsequence formed by
these hook deletions is, very crudely, $\lfloor a/3\rfloor$ at most, and the total number of
those subsequences is at most some $S_2(a)=O(1)$. So $g^{\la}(\boldsymbol{\alpha})$,
the overall number of ways to empty $\la$, is bounded by  the number of $\lfloor N/3
\rfloor$-long $\{0,1\}$-sequences, with at most $\lfloor a/3\rfloor$ $1$'s, corresponding to
the deletions of rim hooks from the bottom subdiagram $\nu$, multiplied by $S_1(a)S_2(a)$,
whence 
\[
g^{\la}(\boldsymbol{\alpha})=O\left(\binom{\lfloor N/3\rfloor}{\lfloor a/3\rfloor}\right)
=O\left(N^{\lfloor a/3\rfloor}\right).
\]
Consequently in the sum $\Sigma_2$, i. e. for $\la_1=N-a$ with $a\le 6$,
\[
\frac{\chi^{\la}({\cal C}_J)\chi^{\la}({\cal C}_2)}{f^{\la}}=
O\left(\frac{N^{\lfloor a/2\rfloor +\lfloor a/3\rfloor}}{N^a}\right).
\]
Since
\[
\min_{a\in [1,6]} \left[a - \left(\lfloor a/2\rfloor +\lfloor a/3\rfloor\right)\right]=1,
\]
and $|\Lambda|=S_3(a)$ is fixed, we obtain then that 
\begin{equation}\label{Sigma2=O}
\Sigma_2=O(N^{-2}).
\end{equation}
Combining \eqref{Sigma1=o} and \eqref{Sigma2=O}, we obtain
\[
\|P_{\gamma}-U\|^2 =O(N^{-2}).
\]
The proof of Theorem \ref{conj2} is complete.
\end{proof}

\section{Number of vertices and Euler characteristic}\label{s:num-vert}

The next claim is directly implied by Theorem \ref{conj2}.

\begin{Theorem}\label{vertexnum} Let $V_N$ denote the number of vertices on the surface
formed by randomly gluing polygons with the sides numbers from $J$, $\min J\ge 3$, such that the counts
$n_j$ of polygons with $j$ sides satisfy $\sum_{j\in J} jn_j = N$. Let $C_N^e$, ($C_N^o$ resp.)
denote
the total number of cycles of the permutation chosen uniformly at random from all even (all
odd permutations resp.) of $[N]$. If  $\alpha$ and $\beta$ are of the same parity (the opposite
parity resp.), then
\[
\|P_{V_N}-P_{C_N^e}\|=O(N^{-1}),\quad (\|P_{V_N}-P_{C_N^o}\|=O(N^{-1}) \text{ resp.}),
\]
uniformly for all admissible $\{n_j\}_{j\in J}$.
\end{Theorem}
\noindent {\bf Note.\/} Fleming and Pippenger \cite{FlemingPippenger} used Gamburd's 
Theorem \ref{Gam1} to  evaluate the $\ell$-th central moment of $V_N$ for $J=\{k\}$, $k\ge 3$,
within an additive error term $O(N^{-1/12}\ln ^{\ell} N)$, ($\ell\ge 1$),  for $J=\{k\}$, $k\ge 3$,
and $N$ divisible by $2\,\text{lcm}\{2,k\}$. With Theorem \ref{vertexnum} at hand, the estimates
in \cite{FlemingPippenger} can be extended to all $N$ divisible by $\text{lcm}\{2,k\}$ with 
a smaller error term $O(N^{-1}\ln^{\ell}N)$.\\

Let us have a look at  $P_{C_N^{e,o}}$. Let $s(N,\ell)$ be the {\it signless\/} Stirling number of
first kind, i. e. the number of permutations of $[N]$ with $\ell$ cycles. Then
\begin{equation}\label{PCNe}
\pr(C_N^e=\ell)= \frac{2s(N,\ell)}{N!},\text{ if }N-\ell\text{ even}; \text{ else }\pr(C_N^e=\ell)=0,
\end{equation}
and
\begin{equation}\label{PCNo}
\pr(C_N^o=\ell)= \frac{2s(N,\ell)}{N!},\text{ if }N-\ell\text{ odd};\text{ else }\pr(C_N^o=\ell)=0,
\end{equation}
see, for instance, Sachkov and Vatutin \cite{SachkovVatutin}. (The equation \eqref{PCNe} is
implicit in Fleming and Pippenger \cite{FlemingPippenger}, Equation (2.2).) Thus the ranges of $C_N^e$ and $C_N^o$ interlace each other. Now $\{s(N,\ell)/N!\}_{\ell\le N}$ is
distribution of $C_N$, the number of cycles in the random permutation of $[N]$, 
and it is well known that $C_N$ is asymptotically normal with mean and variance given by
\begin{align*}
\ex[C_N]&=\sum_{j=1}^N \frac{1}{j}=\ln N+O(1), \\
\var(C_N)&=\sum_{j=1}^N \frac{1}{j}\left(1-\frac{1}{j}\right)=\ln N+O(1).
\end{align*}
The standard proof is based on the observation that $C_N$ has the same distribution as $\sum_{j=1}Y_j$, where $Y_j\in \{0,1\}$ are independent with $\pr(Y_j=1)=1/j$. In fact, Kolchin \cite{Kolchin} had proved a local limit theorem (LLT) for $C_N$, which implies the integral asymptotic normality of $C_N$:
 \begin{equation}\label{CNlocal}
\pr(C_N=\ell) =\frac{(1+o(1))\exp\left(-\frac{(\ell-E[C_N])^2}{2\var(C_N)}\right)}{\sqrt{2\pi\var(C_N)}} , 
\end{equation}
uniformly for $\ell$ such that 
\begin{equation}\label{+/-an12}
\frac{\ell - \ex[C_N]}{\sqrt{\var(C_N)}}\in [-a , 
a],\quad a>0 \text{ fixed}.
\end{equation}
Note that the probability generating function of $C_N$ has only real roots $0,-1,\dots, -(n-1)$, so
that, by Menon's theorem \cite{Menon}, the distribution of $C_N$ is log-concave. Using Canfield's
quantified version of Bender's LLT for log-concave distributions (\cite{Bender}, \cite{Canfield}),
one can show that in \eqref{CNlocal} $o(1)=O(\var(C_N)^{-1/4})$. Applying a  LLT proved
recently in Lebowitz et all \cite{Lebowitz} this bound can be further improved to
$O(\var(C_N)^{-1/2})$.  Combining \eqref{CNlocal} with \eqref{PCNe}-\eqref{PCNo}, we obtain: uniformly 
for $\ell$ satisfying \eqref{+/-an12},
\begin{equation}\label{pr(CNeo=ell)=}
\pr(C_N^{e,o}=\ell) =\frac{\bigl(2+O(\var(C_N)^{-1/2})\bigr)\exp\left(-\frac{(\ell-E[C_N])^2}{2\var(C_N)}\right)}{\sqrt{2\pi\var(C_N)}} ;
\end{equation}
here $N-\ell$ is even (odd resp.) for $C_N^e$ ($C_N^o$ resp.) The equation \eqref{pr(CNeo=ell)=}
and Theorem \ref{vertexnum} taken together 
imply a strong {\it local\/} limit theorem for $V_N$.
\begin{Theorem}\label{LocVN} Uniformly for all admissible $\ell$, meeting \eqref{+/-an12},
\begin{equation}\label{pr(VN=ell)=}
\pr(V_N=\ell)=\frac{\bigl(2+O(\var(C_N)^{-1/2}))\bigr)\exp\bigl(-\frac{(\ell-\ex[C_N])^2}{2\var(C_N)}\bigr)}{\sqrt{2\pi \var(C_N)}};
\end{equation}
admissibility means that $N-\ell$ is even (odd resp.) when $\alpha$ and $\beta$ are
of the same parity (the opposite parity resp.). Consequently $V_N$ is asymptotically normal
with mean and variance $\ln N$ both, $V_N\sim {\cal N}(\ln N,\ln N)$ in short.
\end{Theorem}
\noindent {\bf Note.\/} Gamburd used his Theorem \ref{Gam1} to prove that $V_N$ is asymptotic in distribution (i. e. integrally) to ${\cal N}(\ln N,\ln N)$ for $J=\{k\}$ and $N$ divisible by $2\,\text{lcm}\{2,k\}$.\\

Since the surface has $V_N$ vertices, $N/2$ edges and $n=\sum_j n_j$ faces,
its Euler characteristic $\chi_N$ is
\begin{equation*}
\chi_N =V_N-N/2+n.
\end{equation*}
Using Theorem \ref{LocVN}, we obtain then
\begin{Corollary}\label{pr(chi=)} \begin{equation}\label{pr(chi=)sim}
\pr\left(\chi_N=-N/2+n +\ell\right)=\frac{\bigl(2+O(\var(C_N)^{-1/2})\bigr)\exp\bigl(-\tfrac{(\ell-\ex[C_N])^2}{2\var(C_N)}\bigr)}{\sqrt{2\pi \var(C_N)}},
\end{equation}
uniformly for all admissible $\ell$, satisfying \eqref{+/-an12}.

\end{Corollary}
\noindent {\bf Note.\/} In effect, the equation \eqref{pr(chi=)sim} gives an asymptotic formula for the {\it fraction\/}
of surfaces with a given value of the Euler characteristic in the case when the absolute-value difference between the number of vertices and $\ln N$ is of order $O\bigl((\ln N)^{1/2}\bigr)$. 

\section{Number of components}\label{s:num-comp} 
Let $X_N$ denote the total number of components of the random surface. 
\begin{Theorem}\label{connect} 
\[
\pr(X_N=1)=1-O(N^{-1}).
\]
\end{Theorem}
\noindent {\bf Notes.\/}  {\bf (1)} This estimate is qualitatively best  in general, since Pippenger and Schleich 
\cite{PippengerSchleich} proved that $\pr(X_N=1)= 1- 5/(6N)+O(N^{-2})$ for $J=\{3\}$. {\bf (2)\/} $X_N$ can be
viewed as the number of components in a random {\it multigraph\/} $MG$ on $n=\sum_jn_j$ vertices,
with the given vertex-degree sequence, such that $n_j$ vertices have degree $j$. 
(Each of the vertices $j$ is represented by a set $S_j$ of cardinality $j$, and two vertices
$j$ and $j'$ are joined by an edge iff in the uniformly random matching $M$ on $S=\uplus_jS_j$ there are points $u\in S_j$ and $u'\in S_{j^\prime}$ such that $(u,u^\prime)\in M$. This model was introduced by Bollob\'as \cite{Bollobas}.) The theorem
\ref{connect} asserts that $MG$ is connected  with probability $1-O(N^{-1})$, uniformly over {\it all\/} degree sequences bounded by $3$ from below. For the maximum degree $\le n^{0.02}$ this
claim is implicit in \L uczak \cite{Luczak}, its focus being on graphs, rather than multigraphs;
see also an earlier result by Wormald \cite{Wormald} for the bounded maxdegree case. 
\begin{proof}
If $X_N>1$ then there exists a partition of the $n=\sum_{j\in J} n_j$ cycles into two groups such that no
two sides of a pair of cycles belonging to different groups are glued together; call it  ``no-match''
condition. A generic 
partition into two groups of cycles is given by the two sets, $\{n_j^\prime\}_{j\in J}$ and $\{n_j^{\prime\prime}\}_{j\in J}$, such that $n_j^\prime + n_j^{\prime\prime}=n_j$, $j\in J$. Introduce
$N^\prime=\sum_j jn_j^\prime$, $N^{\prime\prime}=\sum_j jn_j^{\prime\prime}$; so $N=N^\prime
+N^{\prime\prime}$. For an admissible partition, both  $N^\prime$ and $N^{\prime\prime}$ must
be even. Consequently $N^\prime,\,N^{\prime\prime}\ge m$, where $m=\min J$ if $\min J$ is even,  and $m=2\min J$ otherwise. The probability of no-match is
\[
P(N^\prime,N^{\prime\prime}):=\frac{(N^\prime-1)!!(N^{\prime\prime}-1)!!}{(N-1)!!}.
\]
Using Stirling formula, we obtain: uniformly for $N^\prime\ge m$, $N^{\prime\prime}\ge
m$,
\begin{equation}\label{P(N'N'')appr}
P(N^\prime,N^{\prime\prime})=O(P^*(N^\prime,N^{\prime\prime})),\quad 
P^*(N^\prime,N^{\prime\prime})=\frac{(N^\prime)^{N^\prime/2}(N^{\prime\prime})
^{N^{\prime\prime}/2}}{N^{N/2}}.
\end{equation}
Furthermore, the total number of $\{n_j^\prime,n_j^{\prime\prime}\}_{j\in J}$ with parameters
$N^\prime$, $N^{\prime\prime}$ is given by
\begin{align}
Q(N^\prime,N^{\prime\prime})&=\sum_{ n_j^\prime+n_j^{\prime\prime}=n_j\atop \sum_jjn_j^\prime=N^\prime;\,\,\sum_j j n_j^{\prime\prime}=N^{\prime\prime}}\prod_{j\in J}\frac{n_j!}{n_j^{\prime}!\,\,n_j^{\prime\prime} !}\notag\\
&=\bigl[x_1^{N^\prime}\,x_2^{N^{\prime\prime}}\bigr]\prod_{j\in J}\sum_{n_j^\prime+n_j^{\prime\prime}=n_j}\frac{n_j!}{n_j^{\prime}!\,\,n_j^{\prime\prime} !}\,x_1^{jn_j^\prime}x_2^{jn_j^{\prime\prime}}\notag\\
&=\bigl[x_1^{N^\prime}\,x_2^{N^{\prime\prime}}\bigr]\prod_{j\in J}(x_1^j+x_2^j)^{n_j}.
\label{QN',N''=coef}
\end{align}
Now
\[
\pr(X_N>1)\le \sum_{N^\prime,\,N^{\prime\prime}\ge 2\atop N^\prime +N^{\prime\prime}=N}
P(N^\prime,N^{\prime\prime})Q(N^\prime,N^{\prime\prime}),
\]
so, by \eqref{P(N'N'')appr}, we need to bound  $P^*(N^\prime,N^{\prime\prime})Q(N^\prime,N^{\prime\prime})$ for the generic $N^\prime$, $N^{\prime\prime}$. By symmetry, it suffices to consider $N^\prime\le N^{\prime\prime}$. By \eqref{QN',N''=coef}
and $N^\prime+N^{\prime\prime}=\sum_j jn_j$, we have: for all $x_1>0$, $x_2>0$,
\begin{equation}\label{QN'N''<}
\begin{aligned}
Q(N^\prime,N^{\prime\prime})&\le x_1^{-N^\prime}x_2^{-N^{\prime\prime}}
\prod_{j\in J}(x_1^j+x_2^j)^{n_j}\\
&\le y^{-N^{\prime}}\prod_{j\in J}(y^j+1)^{n_j}=\exp(H(y,N^\prime));\\
H(y,N^{\prime})&:=\sum_j n_j\ln (y^j+1) - N^\prime \ln y,\quad y:=\frac{x_1}{x_2}.
\end{aligned}
\end{equation}
The best value of $y$ minimizes $H(y,N^\prime)$, and so it is a root of $H_y(y,N^\prime)=0$,
which is equivalent to
\begin{equation}\label{yrootof}
\sum_j n_j\frac{jy^j}{y^j+1}=N^\prime.
\end{equation}
The LHS strictly increases with $y$, and equals $0$ at $y=0$ and $\sum_j jn_j/2=N/2\ge N^{\prime}$. So $H(y,N^\prime)$ does attain its minimum at a unique point $y(N^\prime)\in (0, 1]$ for
all $N^\prime\le N/2$. $y(N^\prime)$ is strictly increasing with $N^\prime$, and---considered
as a function of the continuous parameter $N^\prime$---$y(N^\prime)$ is {\it continuously
differentiable\/} for $N^\prime>0$, as
\[
\frac{d}{dy}\left(\sum_j n_j\frac{jy^j}{y^j+1}\right) >0,\quad \forall\, y>0.
\]
Thus $Q(N^\prime,N^{\prime\prime})\le \exp\bigl(H(y(N^\prime),N^\prime)\bigr)$, and so
\begin{equation}\label{P*Q=OexpH}
\begin{aligned}
P^*(N^\prime,N^{\prime\prime})Q(N^\prime,N^{\prime\prime})&=O\bigl(\exp({\cal H}(y(N^\prime),N^\prime))\bigr),\\
{\cal H}(y(N^\prime),N^\prime)):&=H(y(N^\prime),N^\prime)\\
&\quad+(N^\prime/2)\ln N^\prime+
(N^{\prime\prime}/2)\ln N^{\prime\prime}-(N/2) \ln N.
\end{aligned}
\end{equation}
We want to show that ${\cal H}(y(N^\prime),N^\prime))$ is strictly decreasing with $N^\prime$.
Since $H_y(y,N^\prime)\bigr|_{y=y(N^\prime)}=0$, and $N^{\prime\prime}=N-N^\prime$, we have
\begin{align*}
\frac{d}{dN^\prime}\,{\cal H}(y(N^\prime),N^\prime))&=\frac{\partial}{\partial N^\prime}\,{\cal H}(y,N^\prime))
\bigr|_{y=y(N^\prime)}\\
&=-\ln y(N^\prime) +(1/2)\ln N^\prime-(1/2)\ln N^{\prime\prime}\\
&=\ln\left(\sqrt{\frac{N^\prime}{N^{\prime\prime}}}\cdot \frac{1}{y(N^\prime)}\right).
\end{align*}
Therefore we need to show that $y(N^\prime)> y_1=y_1(N^\prime):=\sqrt{\tfrac{N^\prime}{N^{\prime\prime}}}$ for $N^\prime<N/2$, or equivalently by \eqref{yrootof}, that 
\[
\sum_j n_j\,\frac{jy_1^j}{y_1^j+1}<N^\prime.
\]
By convexity of $z/(1+z)$ for $z\ge 0$,  
\begin{align*}
\sum_j n_j\,\frac{jy_1^j}{y_1^j+1}&\le N\frac{\sum_j y_1^j(j n_j)/N }{\sum_j y_1^j(j n_j)/N +1}\\
&\le N\frac{y_1^3\sum_j (j n_j)/N}{y_1^3\sum_j (j n_j)/N+1}\\
&=N\,\frac{y_1^3}{y_1^3+1}.
\end{align*}
and
\begin{align*}
\sum_j n_j\,\frac{jy_1^j}{y_1^j+1}-N^\prime&\le (N^\prime+N^{\prime\prime})\,\frac{y_1^3}{y_1^3+1}-N^\prime\\
&=N^{\prime\prime}\left[(y_1^2+1)\,\frac{y_1^3}{y_1^3+1} - y_1^2\right]\\
&=-N^{\prime\prime}\,\frac{y_1^2(1-y_1)}{y_1^3+1}<0,
\end{align*}
for $N^\prime< N/2$. Thus indeed $y(N^\prime)>y_1$, whence  ${\cal H}(y(N^\prime),N^\prime))$ is strictly decreasing for $N^\prime\in (0,N/2]$.

Consider $N^\prime\in [\nu, N/2]$, $\nu=\lfloor 6\ln N\rfloor$, so that $N^{\prime\prime}=N-\nu$.
Then, with $y_1:=y_1(\nu)$, $n=\sum_j n_j$,
\begin{align*}
{\cal H}(y(N^\prime),N^\prime))&\le {\cal H}(y(\nu),\nu)\le {\cal H}(y_1,\nu))\\
&=\sum_j n_j\ln(1 + y_1^j)-(\nu/2) \ln (\nu/N^{\prime\prime})\\
&\,\,\,\,\,+(\nu/2)\ln \nu+(N^{\prime\prime}/2)\ln N^{\prime\prime} -(N/2)\ln N\\
&\le \sum_j n_j y_1^j -(N/2)\ln\bigl(N/N^{\prime\prime}\bigr)\\
&\le n\left(\frac{6\ln N}{N^{\prime\prime}}\right)^{3/2}-\frac{\nu}{2}=O\bigl(\sqrt{N^{-1}\ln ^3N}\,\bigr)-
\frac{\nu}{2}\\
&\le -2 \ln N.
\end{align*}
From this bound and \eqref{P*Q=OexpH} it follows then that
\begin{equation}\label{N'>6logN}
\sum_{N^\prime+N^{\prime\prime}=N\atop 6\log N\le  N^\prime\le N^{\prime\prime}}\!\!\!\!\!\!\!
P^*(N^\prime,N^{\prime\prime})Q(N^\prime,N^{\prime\prime})=O\bigl(N\exp(-2\log N)\bigr)=O(N^{-1}).
\end{equation}
It remains to consider $m\le N^\prime \le 6\ln N$. We will use the bound \eqref{QN'N''<} again, but
this time we are content with a suboptimal $\hat y=\hat y(N^\prime):=(N^\prime/N^{\prime\prime})^{1/j_1}$, $j_1:=\min J\ge 3$. Using the resulting bound for $P^*(N^\prime,N^{\prime\prime})Q(N^\prime,N^{\prime\prime})$, i. e. \eqref{P*Q=OexpH} with $\hat y(N^\prime)$ instead of $y(N^\prime)$, and also 
\[
\sum_j n_j\ln (1+\hat{y}^j)\le \sum_j \hat{y}^j n_j\le \frac{N^\prime}{N^{\prime\prime}}\sum_j n_j=
\frac{n N^\prime}{N^{\prime\prime}},
\]
we obtain
\begin{align*}
P^*(N^\prime,N^{\prime\prime})Q(N^\prime,N^{\prime\prime})&=O\bigl(\exp(\hat{\cal H}(N^\prime))
\bigr),\\
\hat{\cal H}(N^\prime)&:=\frac{n N^\prime}{N^{\prime\prime}}
+ \left(\frac{1}{2}-\frac{1}{j_1}\right)N^\prime\ln N^\prime\\
&+\left(\frac{N^\prime}{j_1}
+\frac{N^{\prime\prime}}{2}\right)\ln N^{\prime\prime}-\frac{N}{2}\ln N.
\end{align*}
Considering  $N^\prime$ as a continuously varying parameter, 
\begin{align*}
\frac{d \hat{\cal H}}{d N^\prime}&=\frac{nN}{(N^{\prime\prime})^2}-\left(\frac{1}{2}-\frac{1}{j_1}\right)
\ln\frac{N^{\prime\prime}}{N^\prime}
 -\frac{1}{j_1}-\frac{N^\prime}{j_1N^{\prime\prime}}\\
&\le -0.5\left(\frac{1}{2}-\frac{1}{j_1}\right)\ln N, 
\end{align*}
uniformly for $m\le N^\prime \le 6\ln N$. Therefore 
\begin{align*}
\sum_{m\le N^\prime\le 6\ln N}\exp(\hat{\cal H}(N^\prime))&\le \exp(\hat{\cal H}(m))
\sum_{s\ge 0}\left[\exp\bigl(-0.5(1/2-1/j_1)\ln N\bigr)\right]^s\\
&\le 2\exp(\hat{\cal H}(m))=O\bigl(N^{-m(1/2-1/j_1)}\bigr).
\end{align*}
If $j_1\ge 3$ is even then $m=j_1\ge 4$, and if $j_1$ is odd then $m=2j_1\ge 6$; so
$m(1/2-1/j_1)\ge 1$. Consequently
\begin{equation}\label{N'<6logN}
\sum_{m\le N^\prime\le 6\ln N}\exp(\hat{\cal H}(N^\prime))=O(N^{-1}).
\end{equation}
Combining \eqref{N'<6logN} with \eqref{N'>6logN}, we complete the proof of Theorem \ref{connect}.
\end{proof}
On the event $\{X_N=1\}$, the genus $g_N$ is given by $g_N=1-\chi_N/2$. Thus $g_N$ is defined with
probability $1-O(N^{-1})$, and so by Corollary \ref{pr(chi=)} we have
\begin{Corollary}\label{pr(g=)} For all admissible $\ell$, and $N$ large enough,
\begin{equation}\label{pr(g=)sim}
\pr\left(g_N=1+\frac{N}{4}-\frac{n}{2}-\frac{\ell}{2}\right)=\frac{\bigl(2+O(\var(C_N)^{1/2})\bigr)\exp\bigl(-\tfrac{(\ell-\ex[C_N])^2}{2\var(C_N)}\bigr)}{\sqrt{2\pi \var(C_N)}},
\end{equation}
($n=\sum_jn_j$), uniformly for all admissible $\ell$, satisfying 
$\tfrac{\ell - \ex[C_N]}{\sqrt{\var(C_N)}}\in [-a, a]$.
\end{Corollary}
\noindent In particular, for $J=\{k\}$, $k\in [3,N]$,
\[
1+\frac{N}{4}-\frac{n}{2}= 1 +\frac{N(k-2)}{4k};
\]
so that for a single disc with $N$ sides, ($N$ even), the genus $g_N$ is asymptotic, integrally {\it and\/} 
locally, to $N/4 -\tfrac{1}{2}\,{\cal N}(\ln N,\ln N)$. We proved this result in Chmutov and
Pittel \cite{ChmutovPittel} by using the Harer-Zagier \cite{HarerZagier} formula for the generating function of chord diagrams enumerated by the genus of the attendant surface. That study
was prompted by an earlier result of Linial and Nowik
\cite{LinialNowik}, who proved, using the H-Z formula,  that $\ex[g_N]= N/4-0.5\ln N+O(1)$.
(They also proved that $\ex[g_N]=N/2 - \Theta(\ln N)$ for a different random surface induced by an
oriented chord diagram, for which a counterpart of the H-Z formula is unknown.)

\end{document}